\input amstex
\documentstyle{amsppt}
\magnification=\magstep1 \baselineskip=12pt \hsize=6truein
\vsize=8truein

\topmatter

\title
Connected Sums of Closed Riemannian Manifolds and Fourth Order
Conformal Invariants
\endtitle
\author David Raske \\
\endauthor

\leftheadtext{Connected sum} \rightheadtext{David Raske}

\email david.t.raske{\@}gmail.edu \endemail

\abstract In this note we take some initial steps in the
investigation of a fourth order analogue of the Yamabe problem in
conformal geometry. The Paneitz constants and the Paneitz
invariants considered are believed to be very helpful to
understand the topology of the underlined manifolds. We calculate
how those quantities change, analogous to how the Yamabe constants
and the Yamabe invariants do, under the connected sum operations.
\endabstract

\endtopmatter

\def\p {\frac{n}{n-4}}
\def \q {\frac{2n}{n-4}}

\def \CM {{C_+}^{\! \! \! \infty}(M)}
\def \U { u_{{\delta}_1}}
\def \V { u_{{\delta}_2}}
\def \Y { ((\Delta F)^2 + a_n R |\nabla F|^2 - \frac{4}{n-2}
\text{Ric} (\nabla F,\nabla F) + Q F^2) \quad dv}
\def \Z { ((\Delta F_l)^2 + a_n R |\nabla F_l|^2 - \frac{4}{n-2}
\text{Ric} (\nabla F_l,\nabla F_l) + Q F^2) \quad dv}
\def \W{ ((\Delta f)^2 + a_n R |\nabla f|^2 - \frac{4}{n-2}
\text{Ric} (\nabla f,\nabla f) + Q f^2) \quad dv}
\def \Q{ ((\Delta f_l)^2 + a_n R |\nabla f_l|^2 - \frac{4}{n-2}
\text{Ric} (\nabla f_l,\nabla f_l) + Q {f_l}^2) \quad dv}

\document
\vskip 0.1in \noindent{\bf 1. Introduction}\vskip 0.1in

Let $(M,g)$ be a connected compact Riemannian manifold without
boundary of dimension $n \geq 5$. Let
$$
Q[g] = - \frac {n-4}{4(n-1)} \Delta R + \frac {(n-4)(n^3 - 4n^2 +
16n - 16)}{16(n-1)^2(n-2)^2}R^2 - \frac {2(n-4)}{(n-2)^2}|Ric|^2
\tag 1.1
$$
be the so-called $Q$-curvature, where $R$ is the scalar curvature,
$Ric$ is the Ricci curvature. And let
$$
P[g] = (-\Delta)^2 - \text{div}_g((\frac {(n-2)^2 +
4}{2(n-1)(n-2)} Rg - \frac 4{n-2}Ric_g)d) + Q[g] \tag 1.2
$$
be the so-called the Paneitz-Branson operator. It is known that
$$
P[g] u = Q[g_u] u^{\frac {n+4}{n-4}} \tag 1.3
$$
which is called the Paneitz-Branson equation, where $g_u =
u^{\frac 4{n-4}}g$ (cf. [P] [Br] [XY] [DHL] [DMA] ). We consider
the equation (1.3) as a fourth order analogue of the well-known
scalar curvature equation
$$
L[g] v= R[g_v] v^{\frac {n+2}{n-2}}, \tag 1.4
$$
where
$$
L[g] = - \frac {4(n-1)}{n-2}\Delta + R \tag 1.5
$$
is the so-called conformal Laplacian and $g_v = v^{\frac
4{n-2}}g$. The well-known Yamabe problem in conformal geometry is
to find a metric, in a given class of conformal metrics, which is
of constant scalar curvature, i.e. to solve
$$
L[g] v =Y v^{\frac {n+2}{n-2}}
$$
on a given manifold $(M, g)$ for some positive function $v$ and a
constant $Y$. The affirmative resolution to the Yamabe problem was
given in [Sc] after other notable works [Ya] [Tr] [Au]. In fact,
it was proven that there exists a so-called Yamabe metric $g_v$ in
the class $[g]$ which is a minimizer for the so-called Yamabe
functional
$$
Y(v) = \frac {\int_M (vL[g]v)dv_g}{(\int_M v^\frac {2n}{n-2}
dv_g)^\frac {n-2}n}.
$$

In chapter one we investigate a fourth order analogue of the
Yamabe problem. Let $\CM$ be the space of smooth non-negative
functions on $M$. Similar to the Yamabe problem, we define the
Paneitz functional
$$
\wp_g(u) = \frac { \int_M (uP[g]u) dv_g}{(\int_M u^\q dv_g)^{\frac
{n-4}n}} \tag 1.6
$$
for $u \in \CM$ and the {\it Paneitz constant} associated with
$(M, [g])$
$$
\lambda(M, [g]) =\inf_{u\in  \CM} \wp (u). \tag 1.7
$$
It is clear that $\lambda(M, [g])$ is a conformal invariant of the
conformal class $[g]$ because of the conformally covariant
property of the Paneitz-Branson operator:
$$
P[g_w]u = w^{-\frac {n+4}{n-4}} P[g] (w\cdot u) \tag 1.8
$$
where $g_w = w^\frac 4{n-4} g \in [g]$. To describe the
differential structure of $M$, we define
$$
\lambda (M) = \sup_{[g]} \lambda(M, [g]). \tag 1.9
$$
We will refer to $\lambda (M)$ as the {\it Paneitz Invariant} of
the manifold $M$ as the counter part of Yamabe invariant. In [Gi],
Gil-Medrano studied the Yamabe constant for a connected sum of two
closed manifolds. One interesting consequence of connected sum
results in [Gi] is that every compact manifold without boundary
admits a conformal class of metrics whose Yamabe constant is very
negative. In Section 2 of Chapter One we calculate as Gil-Medrano
did in [Gi] to verify that

\proclaim{Theorem 1.1} Let $(M_1,g_1)$ and $(M_2,g_2)$ be two
compact Riemannian manifolds of dimension $n \geq 5$. Then, for
each $\epsilon > 0$, there is a conformal class $[g]$ of metrics
on $M_1 \# M_2$ such that
$$
\lambda(M_1 \# M_2,[g]) < \min\{\lambda(M_1, [g_1]),\lambda(M_2,
[g_2])\} + \epsilon \tag 1.10
$$
and there exists a conformal class $[h]$ of metrics on $M_1 \#
M_2$ such that
$$
\lambda(M_1 \# M_2, [h]) < 2^{-\frac {n-4}n}(\lambda(M_1, [g_1]) +
\lambda(M_2, [g_2])) + \epsilon. \tag 1.11
$$
\endproclaim

Due to the works of Schoen and Yau [SY] (see also [GL]), one knows
that there is some topological constraint for a manifold to
possess a metric of positive Yamabe constant. Therefore it is
interesting to see how the Yamabe invariant is effected  by
connected sum. It was proven in [Ko] [SY] [GL] that the Yamabe
invariant of connected sum of two manifolds with positive Yamabe
invariants is still positive. More precisely, Kobayashi in [Ko]
showed that the Yamabe invariant of connected sum of two manifolds
is greater than or equal to the smaller of the Yamabe invariants
of the two. In Section 3 of Chapter 1 we obtain an analogue for
the Paneitz invariant.

\proclaim{Theorem 1.2} If $M_1$ and $M_2$ are compact manifolds of
dimension $n \geq 5$, then
$$\lambda(M_1 \# M_2) \geq \min \{
\lambda(M_1), \lambda(M_2) \}. \tag 1.12
$$ \endproclaim

The positivity of Paneitz invariant in dimension higher than $4$
should be a topological constraint, as indicated by successful
researches in [CY] (references therein) for fourth order analogue
of how Gaussian curvature influences the geometry of surfaces in
dimension $4$. Another testing ground is to consider closed
locally conformally flat manifolds. Then the recent works in [CHY]
[G] indicate to us that the positivity of fourth order curvature
is indeed very informative about the topology of the underlined
manifolds. We would also like to mention the work by Xu and Yang
in [XY] where they demonstrated that positivity of the
Paneitz-Branson operator is stable under the process of taking
connected sums of two closed Riemannian manifolds.

In Section 1 of Chapter 1 we discuss some preliminary facts about
the Paneitz functional. In Section 2 we calculate and verify
Theorems 1.1. In Section 3 we prove Theorem 1.2.

\vskip 0.1in \noindent{\bf 2. Preliminaries}\vskip 0.1in

Recall that the Yamabe constant of any closed manifold of
dimension greater than $2$ is a finite number and the largest
possible Yamabe constant is realized and only realized by the
Yamabe constant of the standard round sphere in each dimension.
The difficult part is to show that the round sphere is the only
one that has the largest Yamabe constant, which was the last step
in the resolution of Yamabe problem solved by Schoen in [Sc] based
on a positive mass theorem of Schoen and Yau . We observe that, by
(1.3),
$$
\int_M (uP[g]u)dv_g = \int_M u Q[g_u]u^{\frac {n+4}{n-4}}dv_g =
\int_M Q[g_u] u^\frac{2n}{n-4}dv_g = \int_M Q[g_u]dv_{g_u}, \tag
2.1
$$
where $g_u = u^\frac 4{n-4}g \in [g]$. Hence
$$
\aligned \int_M (uP[g]u)dv_g & = \int_M ((\frac {(n-4)(n^3 - 4n^2
+ 16n - 16)}{16(n-1)^2(n-2)^2}R^2 - \frac
{2(n-4)}{(n-2)^2}|Ric|^2)dv)[g_u] \\
& \leq (\frac {(n-4)(n^3 - 4n^2 + 16n -
16)}{16(n-1)^2(n-2)^2}\int_M (R^2)dv)[g_u] \endaligned
$$
When we consider a Yamabe metric $g_u$, i.e.
$$
\frac {\int_M (Rdv)[g_u]}{\text{vol}(M, g_u)^\frac {n-2}n} =
Y\text{vol}(M, g_u)^\frac 2n  \leq n(n-1)\text{vol}(S^n,
g_0)^\frac 2n  , \tag 2.2
$$
we have
$$
\aligned \frac {\int_M (uP[g]u)dv_g}{\text{vol}(M, g_u)^\frac
{n-4}n} & \leq \frac {(n-4)(n^3 - 4n^2 + 16n -
16)}{16(n-1)^2(n-2)^2}Y^2 \text{vol}(M, g_u)^\frac 4n  \\& \leq
\frac {(n-4)(n^3 - 4n^2 + 16n - 16)}{16(n-1)^2(n-2)^2}(n(n-1))^2
\text{vol}(S^n, g_0)^\frac 4n \\& = \frac {\int_{S^n}
(Qdv)[g_0]}{\text{vol}(S^n, g_0)^\frac {n-4}n} = \lambda(S^n,
[g_0]).
\endaligned \tag 2.3
$$
Consequently we obtain

\proclaim{Lemma 2.1} Let $(M^n, g)$ be a closed Riemannian
manifold of dimension great than $4$ with nonnegative Yamabe
constant. Then
$$
\lambda(M^n, [g]) \leq \lambda(S^n, [g_0]) \tag 2.4
$$
and the equality holds if and only if $(M, g)$ is conformally
equivalent to the standard round sphere $(S^n, g_0)$.
\endproclaim

On the other hand, by some choices of testing functions similar to
the ones used to estimate the Yamabe functional, we get

\proclaim{Lemma 2.2} Let $(M^n, g)$ be a closed Riemannian
manifold of dimension great than $4$. Then
$$
-\infty < \lambda(M^n, [g]) \leq \lambda(S^n, [g_0]), \tag 2.5
$$
where $g_0$ is the standard round metric on the sphere $S^n$.
\endproclaim

\demo{Proof} The Paneitz constant is easily seen to be bounded
from the below. Because, by (1.2),
$$
\aligned \int_M (uP[g]u)dv &  = \int_M |\Delta u|^2dv  + a_n\int_M
R|\nabla u|^2dv \\ & \quad - \frac 4{n-4}\int_M \text{Ric}(\nabla
u, \nabla u)dv + \int_M Qu^2dv, \endaligned\tag 2.6
$$
where
$$
a_n = \frac {(n-2)^2 + 4}{2(n-1)(n-2)}.
$$
It suffices to estimate (2.3) for nonnegative functions such that
$$
\int_M u^\frac {2n}{n-4} dv = 1.
$$
Hence, By Holder inequality,
$$
\aligned  \int_M (uP[g]u)dv & \geq \int_M |\Delta u|^2dv - C_1
\int_M |\nabla u|^2 dv - C_2 \int_M u^2 dv \\
& \geq \int_M |\Delta u|^2dv - C_1 \int_M (-\Delta u )u dv - C_2
\int_M u^2dv \\& \geq \frac 12 \int_M |\Delta u|^2dv - \frac 12
C_1^2 \int_M u^2 dv - C_2\int_M u^2dv  \\ & \geq - (\frac 12 C_1^2
+ C_2) (\int_M u^\frac {2n}{n-4}dv)^\frac {n-4}n \text{vol}(M,
g)^\frac 4n \\ & \geq - (\frac 12 C_1^2 + C_2) \text{vol}(M,
g)^\frac 4n.
\endaligned \tag 2.7
$$
for some constants $C_1, C_2 >0$ depending on $(M^n, g)$.

To estimate the upper bound we choose to works in a geodesic
normal coordinate in very small geodesic ball $B_{2\epsilon}
\subset M$ and transplant the rescaled round sphere metric. Let
$B_{2\epsilon}(0) \subset R^n$ and
$$
g_{ij}(x) = \delta_{ij} + O(|x|^2), \forall x\in B_{2\epsilon}(0).
\tag 2.8
$$
Let
$$
u_\epsilon (x) = \left\{ \aligned (\frac {2\epsilon^3}{\epsilon^6
+
|x|^2})^\frac {n-4}2 \quad & \forall x\in B_\epsilon (0) \\
0 \quad\quad\quad \quad  & \forall x\notin
B_{2\epsilon}(0)\endaligned\right. \tag 2.9
$$
be a smooth nonnegative function on $M$. Then it is easily
calculated that
$$
\aligned \int_M (u_\epsilon P[g]u_\epsilon)dv & =
\int_{B_{\epsilon}(0)} |\Delta u_\epsilon|^2dx + o(1) \\ & =
\int_{R^n} |\Delta (\frac {2\epsilon^3}{\epsilon^6 + |x|^2})^\frac
{n-4}2|^2 dx + o(1) \\
& = \int_{R^n} |\Delta (\frac 2{1+|x|^2})^\frac {n-4}2|^2 dx +
o(1) \endaligned\tag 2.10
$$
and
$$
\aligned \int_M u_\epsilon^\frac {2n}{n-4} dv & =
\int_{B_\epsilon(0)}u_\epsilon^\frac {2n}{n-4} dx + o(1) \\ & =
\int_{R^n}(\frac {2\epsilon^3}{\epsilon^6 + |x|^2})^ndx + o(1) \\
& =  \int_{R^n}(\frac {2}{1 + |x|^2})^ndx + o(1). \endaligned\tag
2.11
$$
Therefore
$$
\wp (u_\epsilon) = \frac {\int_M (u_\epsilon
P[g]u_\epsilon)dv}{(\int_M u_\epsilon^\frac {2n}{n-4} dv)^\frac
{n-4}n} = \frac {\int_{R^n} |\Delta s|^2 dx} {(\int_{R^n} s^\frac
{2n}{n-4} dx)^\frac {n-4}n} + o(1), \tag 2.12
$$
where $s = (\frac 2{1+|x|^2})^\frac {n-4}2$. Thus, take $\epsilon
\rightarrow 0$, we arrive at
$$
\lambda(M, [g]) \leq \lambda(S^n, [g_0]). \tag 2.13
$$
\enddemo

One interesting question would be whether $(M, g)$ is conformally
equivalent to $(S^n, g_0)$ when $\lambda(M, [g]) = \lambda(S^n,
[g_0])$ without assuming the Yamabe constant of $(M, g)$ is
nonnegative. In other words one would be interested in searching
for some analogue of a positive mass theorem of Schoen and Yau
here if it make any sense.

\newpage\noindent{\bf 3. Connected Sums and the Paneitz
Constant}\vskip 0.1in

In this section we will calculate the Paneitz functional on a
connected sum of two closed manifolds and verify Theorem 1.1. Let
$(M, g)$ be a closed manifold of dimension higher than $4$. Fix a
point $p\in M$ and let
$$
f_\delta = \left\{\aligned 0 \quad & \forall x \in B_\delta (p) \\
                           1 \quad & \forall x  \in M\setminus
B_{2\delta}(p)
\endaligned\right.\tag 3.1
$$
be a family of smooth functions. We may ask
$$
\left\{\aligned 0 \leq & f_\delta \leq 1 \\
|\nabla f_\delta | & < \frac {C_0}{\delta} \\
| \Delta f_\delta | & < \frac {C_0}{\delta^2}
\endaligned\right. \tag 3.2
$$
for some number $C_0 > 0$. First we calculate

\proclaim{Lemma 3.1} Let $(M, g)$ be a closed manifold of
dimension greater than $4$. Let $u\in \CM$ be given. Then
$u_\delta = f_\delta u \in \CM$ and
$$
\wp_g (u_\delta) = \wp_g(u) + o(1) \tag 3.3
$$
as $\delta \rightarrow 0$
\endproclaim

\demo{Proof} We simply calculate, for a fixed $\delta
> 0$, by (2.6) and (3.2),
$$
\aligned \int_M (u_\delta P[g]u_\delta)dv &  = \int_M |\Delta
u_\delta |^2dv + a_n\int_M R|\nabla u_\delta|^2dv \\ & \quad -
\frac 4{n-4}\int_M \text{Ric}(\nabla u_\delta, \nabla u_\delta)dv
+ \int_M Qu_\delta^2dv \\ & = \int_M (uP[g]u)dv + o(1)
\endaligned\tag 3.4
$$
and
$$
\int_M u_\delta^\frac {2n}{n-4}dv = \int_M u^\frac {2n}{n-4}dv +
o(1), \tag 3.5
$$
as $\delta\rightarrow 0$.
\enddemo

Now let us consider the connected sum of two closed Riemannian
manifolds. Let $(M_1,g_1)$ and $(M_2,g_2)$ be two compact
Riemannian manifolds without boundary of dimension $n \geq 5$. For
$x_1 \in M_1$ and $x_2 \in M_2$, let $B_{\delta_1}(x_1) \subset
M_1$ and $B_{\delta_2} (x_2)\subset M_2$ be geodesic balls
respectively. To make the connected sum one simply to take off the
open balls $B_{\frac 12 \delta_1} (x_1)$ and $B_{\frac 12
\delta_2} (x_2)$ from $M_1$ and $M_2$, identify $\partial B_{\frac
12  \delta_1} (x_1)$ with $\partial B_{\frac 12 \delta_2} (x_2)$.
Hence
$$
\aligned M_1\#M_2 & = \\ \left[(M_1\setminus B_{\frac 12
\delta_1}(x_1)) \bigcup (M_2 \setminus B_{\frac 12 \delta_2}
(x_2))\right] & / \{\partial B_{\frac 12\delta_1}(x_1) \sim
\partial B_{\frac 12\delta_2} (x_2)\}.\endaligned \tag 3.6
$$
We may construct a metric $g$ on the connected sum $M_1 \# M_2$
such that $g$ agrees with $g_1$ on $M_1 \setminus B_{\delta_1}
(x_1)$ and $g_2$ on $M_2\setminus B_{\delta_2} (x_2)$. Notice that
topologically $M_1\#M_2$ does not depend on the value of
$\delta_i$ when they are sufficiently small. Now let us calculate
and estimate the Paneitz functional on the connected sum.

\proclaim{Theorem 3.2} Let $(M_1,g_1)$ and $(M_2,g_2)$ be two
closed Riemannian manifolds of dimension $n\geq 5$. Then for each
$\epsilon > 0$, there is a conformal structure $[g]$ on $M_1 \#
M_2$ such that
$$
\lambda(M_1\#M_2, [g]) < \min\{\lambda(M_1, [g_1]), \lambda(M_2,
[g_2])\} + \epsilon. \tag 3.7
$$
Alternatively, we may find a conformal structure $[g]$ on $M_1 \#
M_2$ such that
$$
\lambda(M,[g]) < \lambda(M_1,[g_1]) + \lambda(M_2,[g_2]) 2^{-\frac
{n-4}n} + \epsilon. \tag 3.8
$$
\endproclaim

\demo{Proof} Let us assume that $\lambda(M_1, [g_1]) \leq
\lambda(M_2, [g_2])$ and $\epsilon > 0$ fixed. By the definition
of the Paneitz constant, we know that there is a real number
$\delta > 0$ and a smooth function $u_\delta \in \CM$ such that
$u_\delta$ vanishes on a geodesic ball $B_\delta (x_1)$ of radius
$\delta$ and centered at $x_1 \in M_1$ and such that
$$
\wp_g(u_\delta) < \lambda(M_1,[g_1]) + \epsilon.
$$
Let $g$ be a metric on $M = M_1 \# M_2$ which agrees with $g_1$,
when restricted to $M_1 \setminus B_\delta (x_1)$. And define the
function $\widetilde{u_\delta}$ on $M_1 \# M_2$ as follows:
$$
\left \{ \aligned \widetilde{u_\delta} & = u_\delta \quad
\text{on} \quad M_1 \setminus B_\delta (x_1) \\
\widetilde{u_\delta} & = 0 \quad \text{elsewhere}.\endaligned
\right.
$$
We then have it that
$$
\wp_g(\widetilde{u_\delta}) = \frac { \int_M ({\Delta
\widetilde{u_\delta}}^2 +  a_n R |\nabla \widetilde{u_\delta}|^2
-\frac{4}{n-2} Ric (\nabla \widetilde{u_\delta},\nabla
\widetilde{u_\delta}) + Q {\widetilde{u_\delta}}^2) dv}{(\int_M
{\widetilde{u_\delta}^\q dv})^{\p}}.
$$
Recalling that $u_\delta$ vanishes on $B_\delta (x_1)$ we see that
$$
\wp_g(\widetilde{u_\delta}) = \wp_{g_1}(u_\delta) <
\lambda(M_1,[g_1]) + \epsilon.
$$
Consequently,
$$
\lambda(M,[g]) < \lambda(M_1,[g_1]) + \epsilon = \min
(\lambda(M_1,[g_1]),\lambda(M_2,[g_2])) + \epsilon.
$$

We will now proceed to prove (3.8). First notice that Lemma 3.1
can be use to say that for any fixed $\epsilon > 0, x_1 \in M_1,
x_2 \in M_2,$ we can find two positive reals ${\delta}_1,
{\delta}_2$ and smooth functions $u_{{\delta}_1}, u_{{\delta}_2},$
where $u_{{\delta}_i} \in C^\infty (M_i),$ with the following
properties:
$$
\left\{\aligned u_{{\delta}_1} & = 0  \quad \text{on
$B_{\delta_1} (x_1)$} \\
\wp_{g_1} (u_{{\delta}_1}) & < \lambda(M_1, [g_1]) + \epsilon_1
\endaligned\right.
$$
and
$$
\left\{\aligned u_{{\delta}_2} & = 0  \quad \text{on $
B_{{\delta}_2} (x_2)$}\\ \wp_{g_2} (u_{{\delta}_2}) & <
\lambda(M_2, [g_2]) + \epsilon_1, \endaligned\right.
$$
where $\epsilon_1 = 2^{-n+4/n} \epsilon$. Also, notice that we can
assume without loss of generality that the $L^\q (M)$ norms of
$u_{{\delta}_1}$ and $u_{{\delta}_2}$ are normalized. Using the
same reasoning as in the proof of (3.7), a metric $g$ on $M_1 \#
M_2$ can be constructed such that $g$ agrees with $g_i$ when
restricted to $M_i\setminus B_{{\delta}_i} (x_i)$. Let us consider
now the function $\widetilde{u}$ on $M = M_1 \# M_2$ given by
$$
\widetilde{u} = \left \{\aligned u_{{\delta}_1} & \quad \text{on}
\quad M_1 \setminus B_{{\delta}_1} (x_1) \\ u_{{\delta}_2} & \quad
\text{on} \quad M_2 \setminus B_{{\delta}_2} (x_1)\\  0 & \quad
\text{elsewhere}
\endaligned \right. \tag 3.9
$$
then
$$
\aligned \wp_g (\widetilde{u}) = \frac { {\int}_{M_1 \setminus
B_{{\delta}_1} (x_1)} ((\Delta \widetilde{u})^2 + a_n R  |\nabla
\widetilde{u}|^2 -\frac{4}{n-4} Ric (\nabla \widetilde{u},\nabla
\widetilde{u}) + Q {\widetilde{u}}^2) dv} {({\int}_{M_1 \setminus
B_ {{\delta}_1} (x_1)} {\widetilde{u}}^\q dv + {\int}_{M_2
\setminus B_{{\delta}_2} (x_2)} {\widetilde{u}}^\q dv)^{\p}} + \\
\frac { {\int}_{M_2 \setminus B_{{\delta}_2} (x_2)} ((\Delta
\widetilde{u})^2 + a_n R  |\nabla \widetilde{u}|^2 -\frac{4}{n-2}
Ric (\nabla \widetilde{u},\nabla \widetilde{u}) + Q
{\widetilde{u}}^2) dv} {({\int}_{M_1 \setminus B_ {{\delta}_1}
(x_1)} {\widetilde{u}}^\q dv + {\int}_{M_2 \setminus
B_{{\delta}_2} (x_2)} {\widetilde{u}}^\q dv)^{\p}} \endaligned
$$
Using (3.9) we then obtain
$$
\aligned \wp_g (\widetilde{u}) = \frac { {\int}_{M_1 \setminus
B_{{\delta}_1} (x_1)} ((\Delta \widetilde{\U})^2 + a_n R  |\nabla
\widetilde{\U}|^2 -\frac{4}{n-2} Ric (\nabla \widetilde{\U},\nabla
\widetilde{\U}) + Q {\widetilde{\U}}^2) dv} {({\int}_{M_1
\setminus B_ {{\delta}_1} (x_1)} {\widetilde{\U}}^\q dv +
{\int}_{M_2 \setminus B_{{\delta}_2} (x_2)} {\widetilde{\V}}^\q
dv)^{\p}} +
\\ \frac { {\int}_{M_2 \setminus B_{{\delta}_2} (x_2)} ((\Delta
\widetilde{\V})^2 + a_n R  |\nabla \widetilde{\V}|^2
-\frac{4}{n-2} Ric (\nabla \widetilde{\V},\nabla \widetilde{\V}) +
Q {\widetilde{\V}}^2) dv} {({\int}_{M_1 \setminus B_ {{\delta}_1}
(x_1)} {\widetilde{\U}}^\q dv + {\int}_{M_2 \setminus
B_{{\delta}_2} (x_2)} {\widetilde{\V}}^\q dv)^{\p}} \endaligned
$$
Now, recalling the above stated properties of $\U$ and $\V$, we
may also assume
$$  {\int}_{M_i \setminus B_{{\delta}_i} (x_i)} {u_{{\delta}_i}}^\q dv
= 1,$$ and $$\aligned \wp_{g_i} (u_{{\delta}_i})& = {\int}_{M_i
\setminus B_{{\delta}_i} (x_i)} ({\Delta
\widetilde{u_{{\delta}_i}}}^2 +  a_n R |\nabla
\widetilde{u_{{\delta}_i}}|^2  -\frac{4}{n-2} \text{Ric} (\nabla
\widetilde{u_{{\delta}_i}},\nabla \widetilde{u_{{\delta}_i}}) + Q
\widetilde{u_{{\delta}_i}}^2) dv \\ & < \lambda(M_i,[g_i]) +
\epsilon_1. \endaligned
$$
Thus
$$
\aligned \lambda(M,[g]) & \leq  \wp_g (\widetilde{u}) \\
& < (\lambda(M_1,[g_1]) + \lambda(M_2,[g_2]) + 2\epsilon_1)
2^{-\frac {n-4}n}
\\ &  = (\lambda(M_1,[g_1]) + \lambda(M_2,[g_2])) 2^{-\frac {n-4}n} +
\epsilon.\endaligned
$$
\enddemo

\vskip 0.1in\noindent{\bf 4. Connected Sums and the Paneitz
Invariants}\vskip 0.1in

Kobayashi in [Ko] showed that the Yamabe invariant of connected
sum of two manifolds is greater than or equal to the smaller of
the Yamabe invariants of the two. The aim of this section is to
generalize this result of Kobayashi to the case of compact
manifolds of dimension $n \geq 5$, and with the Yamabe invariant
$Y(M)$ replaced by it's fourth order analogue the Paneitz
invariant $\lambda(M)$. Namely, we have


\proclaim{Theorem 4.1} If $M_1$ and $M_2$ are closed manifolds of
dimension $n \geq 5$. If $\lambda(M_1) > 0$ and $\lambda(M_2) > 0$
then
$$\lambda(M_1 \# M_2) \geq \min \{
\lambda(M_1), \lambda(M_2) \}. \tag 4.1
$$
\endproclaim

We will basically follow the approach taken by Kobayashi in [Ko].
First we consider the Paneitz invariant on the disjoint union of
compact manifolds. Take two $n$-manifolds with conformal
structures, say $(M_1,[g_1])$ and $(M_2,[g_2])$. We write $(M,[g])
= (M_1,[g_1]) \bigsqcup (M_2,[g_2])$ if $M$ is the disjoint union
of $M_1$ and $M_2$, and $g_i = \{ g|_{M_i} ; g \in [g] \}$ for
$i=1,2$. Let $u$ be a smooth non-negative function on $M$. Since
$M$ is the disjoint union of $M_1$ and $M_2$ it follows that we
can write $u = u_1 + u_2$, where $u_i = 0$ on $M_j$, where $i \neq
j$ and where $u_i$ is a non-negative smooth function on $M_i$. If
we assume that $\lambda(M_i,[g_i]) \geq 0$ for $i = 1,2$, then it
can easily be seen that $$\lambda(M,[g]) = \min \{
\lambda(M_1,[g_1]), \lambda(M_2,[g_2])\}.$$ Due to Lemma 2.2, we
can assume that $\lambda(M_1)$ and $\lambda(M_2)$ are finite; and
we can use the above equation to conclude that
$$\lambda(M) = \min \{
\lambda(M_1) , \lambda(M_2) \}.
$$

Let $M$ be a compact manifold of dimension $n \geq 5$, and $p_1$
and $p_2$ two points of $M$. We take off two small balls around
$p_1$ and $p_2$, and then attach a handle instead, the handle
being topologically the product of a line segment and $S^{n-1}$.
The new manifold obtained in this way will be denoted by
$\overline{M}$. Let $M_1$ and $M_2$ be Riemannian manifolds and
let $M_1 \bigsqcup M_2$ denote the disjoint union of $M_1$ and
$M_2$. If $M = M_1 \bigsqcup M_2$ and $p_1$ and $p_2$ are taken
from $M_1$ and $M_2$ respectively, then $\overline{M} = M_1 \#
M_2$. Therefore we see that in order to prove Theorem 4.1 it
suffices to show
$$
\lambda(\overline{M}) \geq \lambda(M).
$$

\demo{Proof of Theorem 4.1} Let $\epsilon$ be an arbitrary
positive number, which will be fixed throughout. First, we take a
metric $g$ on $M$ such that
$$
\lambda(M,[g]) > \lambda(M) - \epsilon. \tag 4.2
$$
Due to continuity considerations we may assume that $[g]$ is
conformally flat around the points $p_1$ and $p_2$. Then there is
a function $\gamma \in C^\infty (M \setminus \{p_1,p_2\})$ and $g
\in [g]$ such that $\widetilde g = e^\gamma g$ is a complete
metric of $M \setminus \{ p_1,p_2\}$ and that each of the two ends
is isometric to the half infinite cylinder $[0,\infty) \times
S^{n-1}(1)$. For convenience, we write
$$
(M \setminus \{ p_1,p_2\},\widetilde{g}) = [0,\infty) \times
S^{n-1}(1) \bigcup (\widetilde{M},\widetilde{g}) \bigcup
[0,\infty) \times S^{n-1}(1),
$$
where $\widetilde{M}$ is the complement of the two cylinders. We
can glue $(\widetilde{M},\widetilde{g})$ and $[0,l] \times
S^{n-1}(1)$, along their boundaries to get a smooth Riemannian
manifold $(\overline{M},g_l)$, where $\overline{M}$ is as
mentioned in the beginning of the section:
$$
(\overline{M},\overline{g_l}) = (\widetilde{M},\widetilde{g})
\bigcup [0,l] \times S^{n-1}(1). \tag 4.3
$$
We then have
$$
\lambda(\overline{M},[g_l]) = \inf_{f>0} \frac {
\int_{\overline{M}} ((\Delta f)^2 + a_n R |\nabla f|^2
-\frac{4}{n-2} \text{Ric} (\nabla f, \nabla f) + Q f^2)
dv}{(\int_{\overline{M}} f^\q dv )^{\p}},
$$
So, take a positive function $f_l \in C^\infty(\overline{M})$ such
that
$$
\int_{\overline{M}} ((\Delta f_l)^2 + a_n R |\nabla f_l|^2 -
\frac{4}{n-2} \text{Ric} (\nabla f_l,\nabla f_l) + Q f^2) dv  <
\lambda(\overline{M},[g_l]) + \frac{1}{l+1} \tag 4.4
$$
and
$$
\int_{\overline{M}} {f_l}^\q dv =1. \tag 4.5
$$

\proclaim{Lemma 4.2} There is a section, say $\{t_l\} \times
S^{n-1}$, in the cylindrical part of $\overline{M}$ such that
$$\int_{ \{t_l\} \times S^{n-1}} ((\Delta f_l)^2 + a_n R |\nabla f_l|^2
- \frac{4}{n-2} \text{Ric} (\nabla f_l,\nabla f_l) + Q f^2) dv <
\frac{B}{l},$$ where $B$ is a constant independent of $l$.
\endproclaim

\demo{Proof} Using (4.4) we have it that
$$\aligned &\int_{S^{n-1} \times [0,l] } \W  \\
&< \lambda(\overline{M},[g_l]) + \frac{1}{1+l} -
 \int_{\widetilde{M}} ((\Delta f_l)^2 +
 a_n R |\nabla f_l|^2 - \frac{4}{n-2}
 \text{Ric} (\nabla f_l,\nabla f_l) + Q {f_l}^2)
 dv. \endaligned$$
 It follows then that it suffices to demonstrate that
 there exists a constant $D$, independent of $l$,
such that
$$
\int_{\widetilde{M}} ((\Delta f_l)^2 + a_n R |\nabla f_l|^2 -
\frac{4}{n-2} \text{Ric} (\nabla f_l,\nabla f_l) + Q {f_l}^2) dv >
D.
$$
Towards this end, we first notice that we can rewrite (4.3) as
follows:
$$
(\overline{M},\overline{g_l}) = (\widetilde{M_1},\widetilde{g_1})
\bigcup [0,l] \times S^{n-1}(1) \bigcup (\widetilde{M_2},
\widetilde{g_2}),$$ where $(\widetilde{M_i},\widetilde{g_i})$, $i
\in \{1,2\}$, is conformal to $(M_i,g_i) \setminus (B_i(p_i),
\delta)$, where $B_i(p_i)$ is a small ball centered at $p_i$ and
$\delta$ is the Euclidean metric. Now, noting that $a_n R +
\frac{4}{n-4} \text{Ric}$ is a strictly positive operator on the
cylindrical component of $\overline{M}$ and that $Q$ is a strictly
positive function on the cylindrical component, we see that we can
write
$$(\widetilde{M_i},\widetilde{g_i}) = (N_i,h_i) \bigcup
({N_i}',{h_i}')$$ where $({N'}_1,{h_1}') \bigcap ([0,l] \bigcup
S^{n-1}) = S^{n-1} \times \{0\}$; $({N'}_2,{h_2}') \bigcap ([0,l]
\bigcup S^{n-1}) = S^{n-1} \times \{l\}$; ${h_i}'$ is conformally
flat; $a_n R_{{h_i}'} - \frac{4}{n-2} {\text{Ric}}_{{h_i}'}$ is a
positive operator pointwise on ${N_i}'$; and $Q_{{h_i}'}$ is
positive on ${N_i}'$. In geometric terms we can think of
$(N'_i,h'_i)$ as a small part of the necks of the connected sum
$\overline{M}$ adjacent to the cylindrical component. We will now
use this refined decomposition of $\overline{M}$ to decompose
$f_l$; that is, we write $f_{l} = f_{l,1} + f_{c,l} + f_{2,l}$,
where $f_{1,l}$ is supported on $\widetilde{M_1}$;$f_{2,l}$ is
supported on $\widetilde{M_2}$; and $f_{c,l}$ is supported on
${N_1}' \bigcup ([0,l] \times S^{n-1}) \bigcup {N_2}'$.
Furthermore we assume that $f_{1,l}$, $f_{2,l}$, and $f_{c,l}$
vanish smoothly at some nonzero distance away from the boundaries
of their respective supports.  We will now see that the energies
$\int_{\overline{M}} f_{1,l} P_g f_{1,l} dv $,
$\int_{\overline{M}} f_{2,l} P_g f_{2,l} dv $, and
$\int_{\overline{M}} f_{c,l} P_g f_{c,l} dv $ are all bounded
below by a constant independent of $l$. First notice that $f_{c,l}
P_g f_{c,l} \geq 0$ on $\overline{M}$, and hence the last integral
listed above is nonzero. Now, notice that due to our assumption
that $f_{i,l}$, $i \in \{1,2\}$, vanish near the boundaries of
their respective supports, we can extend $f_{i,l}$ to a smooth,
non-negative function $f'_{i,l}$ on $M_i$, by defining $f'_{i,l}$
to be zero on $M_i \setminus \widetilde{M_i}$. Lemma 2.1 then
provides us with the existence of negative constants $D_i$ such
that $\int_{M_i} f_{i,l} P_g f_{i,l} f_{i,l} dv \geq D_i
(\int_{M_i} {f_{i,l}}^\q dv)^\p \geq D_i.$ Since $D_i$ is
determined strictly by the conformal structure of $(M_i,g_i)$, the
above bounds are independent of $l$. Putting these three energy
estimates together we have it that there exists a constant $D$
such that
$$\int_{\widetilde{M}} ((\Delta f_l)^2 + a_n R |\nabla f_l|^2
-\frac{4}{n-2} \text{Ric} (\nabla f_l,\nabla f_l) + Q {f_l}^2) dv
> D.$$ As a consequence we have it that there is a $t_l \in [0,l]$
such that
$$\aligned &\int_{ \{t_l \} \times S^{n-1}} \Q
\\ &< (\lambda(\overline{M},C_l) +\frac{1}{1+l} +
D)/l,
\endaligned$$ which gives us Lemma 4.1 with $B = (\lambda(\overline{M})
+ 1 + B_1)$.

\enddemo

Now we cut off $\overline{M}$ on the section $\{ t_1  \times
S^{n-1} \}$, and attach two half-infinite cylinders to it, so $(M,
\setminus \{p_1, p_2 \}, \overline{g})$ reappears. But this time
we describe it as follows:
$$
(M, \setminus \{p_1, p_2 \}, \overline{g}) = [0,\infty) \times
S^{n-1} (1) \bigcup (\overline{M} - \{ t_1 \} \times S^{n-1} ,
g_l) \bigcup [0,\infty) \times S^{n-1} (1).
$$
We think of the function $f_l$ as defined on $\overline{M} - \{
\{t_l \} \times S^{n-1} \}$, and extend it to the whole space $M -
\{ p_1,p_2 \}$ as follows: Let $F_l$ be Lipschitz function of
$\overline{M} - \{ p_1 , p_2 \}$ such that
$$F_l = f_l \qquad on \quad \overline{M} - \{ t_l \} \times
S^{n-1}$$ and $$F_l (t,x) = \left\{\aligned
(1-t)\widetilde{f_l}(x)
\quad \text{for} \quad (t,x) \in [0,1] \times S^{n-1}; \\
0 \quad \text{for} \quad (t,x) \in [1,\infty] \times S^{n-1},
\endaligned\right.
$$
where $\widetilde{f_l} = f_l|_{\{t_l\}\times S^{n-1}} \in C^\infty
(S^{n-1})$. Now it easy to see from (4.4) and (4.6) that
$$
\int_{ M \setminus \{p_1, p_2\}} \Z < \lambda(\overline{M},[g_l])
+\frac{B}{l},
$$
where $B$ is a constant independent of $l$. Obviously from (4.5)
we get
$$\int_{\overline{M} \setminus \{ p_1, p_2 \}} {F_l}^\q dv > 1.
$$
Therefore, we have
$$
\inf  \frac { \int_{M \setminus \{ p_1,p_2\}} \Y} { (\int_{M
\setminus \{ p_1,p_2\}}  F^\q dv )^{\p}} \leq \lambda
(\overline{M}), \tag 4.9
$$
where the infimum is taken over all nonnegative Lipschitz
functions $F$ with compact support. It follows from the choice of
the metric $\widetilde{g}$ that the left side of (4.9) is equal to
$\lambda (M,[g])$. Since $\epsilon$ can be chosen arbitrarily in
(4.2), we conclude $\lambda(M) \leq \lambda (\overline{M})$, which
completes the proof.
\enddemo

\vskip 0.1in \noindent {\bf References}:

\roster \vskip0.1in \item"{[Au]}" T. Aubin, The scalar curvature,
differential geometry and relativity, (Cahen and Flato, eds.),
Reidel, Dordrecht 1976.

\vskip0.1in \item"{[Br]}" Thomas Branson, Group representations
arising from Lorentz conformal geometry, J. Func. Anal. 74 (1987),
199-291.

\vskip 0.1in \item"{[CHY]}" S.-Y. A. Chang, F. Hang and P. Yang,
On a class of locally conformally flat manifold, Preprint 2003.

\vskip0.1in \item"{[CY]}" S.-Y. A. Chang and P. Yang, Nonlinear
differential equations in conformal geometry, Proceedings of ICM
2002.

\vskip0.1in \item"{[DHL]}" Z. Djadli, E. Hebey and M. Ledoux,
Paneitz type operators and applications, Duke Math. J. 104(2000)
no. 1, 129-169.

\vskip0.1in \item"{[DMA]}" Z. Djadli, A. Malchiodi and M. Ahmedou,
Prescribing a fourth order conformal invariant on the standard
sphere, Part II - blow-up analysis and applications, Preprint
2001.

\vskip 0.1in \item"{[Gi]}" O. Gil-Medrano, Connected sums and the
infimum of the Yamabe functional, Differential geometry, Peñíscola
1985, 160--167, Lecture Notes in Math., 1209, Springer, Berlin,
1986.

\vskip0.1in \item"{[G]}" Maria Gonzalez, Singular sets of a class
of locally conformally flat manifolds, Preprint 2004.

\vskip 0.1in \item"{[GL]}" M. Gromov and H.B. Jr. Lawson,
Classification of simply connected manifolds with positive scalar
curvature, Ann. Math. 111 (1980) 423-434.

\vskip 0.1in \item"{[Ko]}" O. Kobayashi, Scalar curvature of a
metric with unit volume, Math. Ann. 279 (1987) 253-265.

\vskip0.1in \item"{[P]}" S. Paneitz, A quadratic conformally
covariant differential operator for arbitrary pseudo-Riemannian
manifolds, Preprint 1983.

\vskip0.1in \item"{[Sc]}" Rick Schoen, Conformal deformation of a
Riemmanian metric to constant scalar curvature, J. Diff. Geom. 6
(1984), 479-495.

\vskip0.1in \item"{[SY]}" R. Schoen and S.T. Yau, On the structure
of manifolds with positive scalar curvature, Manu. Math. 28 (1979)
159 - 183.

\vskip0.1in \item"{[Tr]}" N. Trudinger, Remarks concerning the
conformal deformation of Riemannian structures on compact
manifolds, Ann. Scuola Norm. Sup. Pisa Cl. Sci (3) 22 (1968),
265-274.

\vskip 0.1in \item"{[XY]}" Xingwang Xu and and Paul Yang,
Positivity of Paneitz operators, Discrete Contin. Dynam. Systems 7
(2001), no. 2, 329--342.

\vskip0.1in \item"{[Ya]}" H. Yamabe, On a deformation of
Riemannian structures on compact manifolds, Osaka J. Math. 12
(1960), 21-37.

\endroster

\enddocument